\documentclass[11pt]{article}
\usepackage{amsfonts,amssymb,amsmath}
\usepackage{amsfonts,amssymb,amsmath,bbm,hyperref}
\hsize \textwidth
\advance \hsize by -\marginparwidth
\evensidemargin \oddsidemargin
\usepackage{amssymb}
\advance\hoffset by 5mm
\newcommand{\pdr}[2]{\frac{\partial{#1}}{\partial{#2}}}
\newcommand{\Rm}{{\mathbb R}}

\newcommand{\expE}{{\mathbb E}}
\def\R{\mathbb{R}}
\newcommand{\eps}{\varepsilon}
\def\epsilon{\varepsilon}

\newcommand{\commentout}[1]{}

\def\e{\varepsilon}

\newcommand{\farc}{\frac}

\newcommand{\br}{\begin{eqnarray}}
\newcommand{\er}{\end{eqnarray}}
\newcommand{\be}{\begin{equation}}
\newcommand{\ee}{\end{equation}}
\newcommand{\baa}{\begin{array}}
\newcommand{\eaa}{\end{array}}
\newcommand{\ba}{\begin{eqnarray}}
\newcommand{\ea}{\end{eqnarray}}
\def\di{\displaystyle}
\def\ln{{\log}}

\numberwithin{equation}{section}

\newtheorem{theorem}{\bf Theorem}[section]

\newtheorem{thm}[theorem]{Theorem}

\newtheorem{lem}[theorem]{Lemma}
\newtheorem{prop}[theorem]{Proposition}

\begin{document}

\title{Convergence to a single wave in the Fisher-KPP equation } 

\author{ James Nolen\thanks{Department of Mathematics, 
Duke University, Durham, NC 27708, USA; nolen@math.duke.edu} \and 
Jean-Michel Roquejoffre\thanks{Institut de Math\'ematiques (UMR CNRS 5219),
Universit\'e Paul Sabatier, 118 route de Narbonne, 31062 Toulouse
cedex, France; jean-michel.roquejoffre@math.univ-toulouse.fr} \and Lenya
Ryzhik\footnote{Department of Mathematics, Stanford University,
Stanford CA, 94305, USA; ryzhik@math.stanford.edu } }

\date{}

%posted to arxiv.org on Monday April 11, 2016.

\maketitle

%%%%%%%%%%%%%%%%%%%%%%%%%%%%%%%%%%%%%%%%%%%%%%%%%%%%%%%%%%%%
%%%%%%%%%%%%%%%%%%%%%%%%%%%%%%%%%%%%%%%%%%%%%%%%%%%%%%%%%%%%
%%%%%%%%%%%%%%%%%%%%%%%%%%%%%%%%%%%%%%%%%%%%%%%%%%%%%%%%%%%%
%%%%%%%%%%%%%%%%%%%%%%%%%%%%%%%%%%%%%%%%%%%%%%%%%%%%%%%%%%%%
\begin{center}
{\it Dedicated to H. Brezis, with admiration and respect}
\end{center}

\begin{abstract}
\noindent We study the large time asymptotics of 
a solution of the  Fisher-KPP reaction-diffusion equation, 
with an initial condition  that is a compact perturbation of a step function. 
A well-known result of Bramson states that, 
in the reference frame moving as $2t - ({3}/{2}) \log t +x_\infty$, 
the solution of the 
equation converges as $t\to+\infty$ to 
a translate of the traveling wave corresponding to the minimal speed~$c_*=2$. 
The constant $x_\infty$ depends on the initial condition $u(0,x)$. The proof is elaborate, and based
on probabilistic arguments. The purpose of this paper is to provide a 
simple proof based on PDE arguments. 
\end{abstract}

\section{Introduction}

We consider  the Fisher-KPP equation:
\begin{equation}
\label{e2.1}
u_t-u_{xx}=u-u^2,~~t>0,~x\in\Rm,
\end{equation}
with an initial condition $u(0,x)=u_0(x)$ which is a compact
perturbation of a step function, in the sense  
that there exist $x_1$ and $x_2$ so that
$u_0(x)=1$ for all $x\le x_1$, and $u_0(x)=0$ for all $x\ge x_2$.
%\begin{equation}
%\label{e2.2}
%u(0,x)=1-H(x)+v_0(x),\  \ \hbox{$v_0$ compactly supported}.
%\end{equation} 
%Here, $H(x)$ is the Heaviside function. 

This equation has a traveling wave solution $u(t,x)=\phi(x-2t)$, moving with the minimal speed~$c_*=2$, connecting the stable 
equilibrium $u\equiv 1$ to the unstable equilibrium~$u\equiv 0$: 
\begin{equation}
\label{e2.4}
\begin{array}{rll}
&-\phi''-2\phi'=\phi-\phi^2,  \\
&\phi(-\infty)=1, \quad \phi(+\infty)=0.
\end{array}
\end{equation}
Each solution $\phi(\xi)$ of 
\eqref{e2.4} is a shift of a fixed profile $\phi_*(\xi)$: 
$
\phi(\xi)=\phi_*(\xi+s),% \quad \forall \;\; \xi \in \Rm,
$
with some fixed $s\in\Rm$. The profile $\phi_*(\xi)$
satisfies the asymptotics
%$$
%\phi_*(\xi)=(\xi+\xi_\infty+k)e^{-(\xi+\xi_\infty)}+O(e^{-(1+\delta)\xi}),~~
%\hbox{ as $\xi\to+\infty$},}
%$$
%for some $\xi_\infty\in\Rm$ depending on $\phi$. In the whole paper, we will fix the 
%translation by selecting the only solution $\phi_*$ of
% \eqref{e2.4} for which $\xi_\infty=0$ in the asymptotics, in other words
 \begin{equation}
\label{e2.20}
\phi_*(\xi)=(\xi+k)e^{-\xi}+O(e^{-(1+\omega_0 )\xi}),~~~\xi\to+\infty,
\end{equation}
with two universal constants $\omega_0 >0$, $k\in\R$.

The large time behaviour of the solutions of this problem has a long history,  starting  with a striking paper of Fisher \cite{Fisher}, which 
identifies the spreading velocity $c_*=2$ via numerical computations and other arguments. 
In the same year, the pioneering KPP paper \cite{KPP} proved 
that the solution of~\eqref{e2.1}, starting from a step function: 
$u_0(x)=1$ for $x\le 0$, $u_0(x)=0$
for $x > 0$, converges to $\phi_*$ in the following sense: there is a   function
\begin{equation}\label{mar2202}
\sigma_\infty(t)=2t+o(t),
\end{equation}
such that 
\[
\di
\lim_{t\to+\infty}u(t,x + \sigma_\infty(t))=\phi_*(x).
\]
Fisher has already made an informal
argument that the $o(t)$ in (\ref{mar2202}) is of the
order $O(\log t)$. 
An important series of papers by   Bramson
 proves the following
\begin{thm}
[\cite{Bramson1}, \cite{Bramson2}]\label{t2.10}
There is  a constant $x_\infty$, depending on $u_0$, such that
\[
\sigma_\infty(t)=2t-\di
\frac32\log t-x_\infty+o(1),\hbox{ as {$t\to+\infty$}}.
\]
\end{thm}
Theorem \ref{t2.10} was proved through elaborate probabilistic arguments.  Bramson also gave necessary and sufficient conditions on the decay of the initial data to zero (as $x \to +\infty$) in order that the solution converges to $\phi_*(x)$ in some moving frame.  Lau~\cite{Lau} also proved those necessary and sufficient conditions (for a more general nonlinear term) using a PDE approach based on the decrease in the number of the intersection points for a pair of solutions of the parabolic Cauchy problem. The asymptotics of $\sigma_\infty(t)$ were not identified by that approach.

A natural question is to prove Theorem \ref{t2.10} with purely PDE arguments. 
In that spirit, a weaker version,
precise up to the~$O(1)$ term, (but valid also for a much more
difficult case of the periodic in space
coefficients), is the main result of~\cite{HNRR1,HNRR2}:
\begin{equation}\label{mar2204}
\hbox{ $\sigma(t)=2t-\di\frac32\log t+O(1)$ as $t\to+\infty$.}
\end{equation} 
Here, we will give a simple and robust proof of Theorem~\ref{t2.10}.
These ideas are further developed 
to study the refined asymptotics of the solutions in~\cite{NRR2}.

The paper is organized as follows. In Section \ref{sec:2}, 
we shortly describe some connections between the Fisher-KPP 
equation \eqref{e2.1} and the branching Brownian motion.
In Section~\ref{sec:3}, we explain, in an informal way, the strategy of 
the proof of the theorem: in a nutshell, the solution
is slaved to the dynamics at $x=O(\sqrt{t})$. In Sections~\ref{sec:4} and \ref{sec:5}, we make the 
arguments of Section~\ref{sec:3} rigorous.

\noindent{\bf Acknowledgment.} JN was  supported by NSF grant DMS-1351653, 
and LR by NSF grant DMS-1311903. JMR  was supported by  the European Union's Seventh 
Framework Programme (FP/2007-2013) / ERC Grant
Agreement n. 321186 - ReaDi - ``Reaction-Diffusion Equations, Propagation and Modelling'',
as well as the ANR project  NONLOCAL ANR-14-CE25-0013.  LR and JMR thank  the Labex CIMI for  a 
PDE-probability  quarter in Toulouse, in Winter 2014, out of which the idea of this
paper grew and which provided a stimulating scientific environment for this project.

\section{Probabilistic links and some related models}\label{sec:2}

The time delay in models of the Fisher-KPP type has been the subject of various 
recent investigations, 
both from the PDE and probabilistic points of view. The Fisher-KPP equation appears 
in the theory of the 
branching Brownian motion (BBM)~\cite{McK} as follows. Consider a BBM
starting at~$x=0$ at time $t=0$, with binary branching at rate $1$. Let $X_1(t),\dots,X_{N_t}(t)$ be
the descendants of the original particle at time $t$, arranged in the
increasing order:
$
X_1(t)\le X_2(t)\le\dots\le X_{N_t}(t).
$
Then, the probability distribution function of the maximum:
$$
v(t,x)={\mathbb P}(X_{N_t}(t)>x), 
$$
satisfies the Fisher-KPP equation
\[
v_t=\di\frac 12 v_{xx}+v-v^2,
\]
with the initial data $v_0(x)=\mathbbm{1}_{x\le 0}$.  Therefore, 
Theorem \ref{t2.10} is about the median location of the maximal particle $X_{N_t}$. 
Building on the work of Lalley and Sellke \cite{LS}, recent 
probabilistic analyses  \cite{ABBS,ABK1,ABK2,BD1,BD2} of this particle 
system have identified a decorated Poisson-type point process which is the limit of 
the particle distribution ``seen from the tip": there is a random variable $Z > 0$ 
such that the point process defined by the shifted particles $\{ X_1(t) - c(t) \, ,\, 
\dots \, , \,X_{N_t}(t)- c(t)\}$, with 
\[
c(t) = 2t - \frac{3}{2} \log t + \log Z,
\]
has a well-defined limit process as $t \to \infty$. Furthermore, 
$Z$ is the limit of the martingale 
\[
Z_t = \sum_{k} (2t - X_k(t))e^{X_k(t) - 2t},
\] 
and 
\[
\hbox{$\phi_*(x) = 1- 
\expE\big[e^{- Ze^{- x}}\big]$ for all $x \in \Rm$.}
\]    

As we have mentioned, the logarithmic term in Theorem \ref{t2.10} arises also
in inhomogeneous variants of this model. For example, consider the Fisher-KPP 
equation in a periodic medium: 
\begin{equation}
\label{e1.300}
u_t-u_{xx}=\mu(x)u-u^2
\end{equation}
where $\mu(x)$ is continuous and 1-periodic in $\R$, such that the 
principal periodic eigenvalue of the operator 
$
-\partial_{xx}-\mu(x)
$
is negative. Then there is a minimal speed $c_* > 0$ such that 
for each $c \geq c_*$, there is a unique  pulsating  front $U_c(t,x)$, up to a time 
shift~\cite{BH,HR}. It was shown in~\cite{HNRR2} that there is~$s_0>0$ such that, if 
$u(t,x)$ solves (\ref{e1.300}) with a nonnegative, nonzero, compactly supported 
initial condition~$u_0(x)$, and~$0<s\leq s_0$, then the $s$-level set $\sigma_s(t)$ of $u(t,x)$ (here, the largest $\sigma>0$ such that $u(t,\sigma) = s$) must satisfy
$$
\sigma_s(t)=c_*t-\frac{3}{2\lambda_*}\log t +O(1),
$$
where $\lambda_* > 0$ is the rate of exponential decay (as $x \to \infty$) of the minimal front $U_{c_*}$, which depends on $\mu(x)$ but not on $s$ or on $u_0$. This implies the convergence of $u(t,x-\sigma_s(t))$ to a closed subset of the family of minimal fronts. It is an open problem to determine whether convergence to a single front holds, not to mention the rate of this convergence.  When $\mu(x) > 0$ everywhere, the solution $u$ of the related model
$$
u_t-u_{xx}=\mu(x)(u-u^2)
$$
may be interpreted in terms of the extremal particle in a BBM with a spatially-varying branching rate \cite{HNRR2}.

Models with temporal variation in the branching process have also been considered. In \cite{FZ}, Fang and Zeitouni studied the extremal particle of such a spatially homogeneous BBM where the branching particles satisfy 
\[
d X(t) = \sqrt{2} \kappa(t/T) \,dB(t)
\]
between branching events, rather than following a standard Brownian motion. 
In terms of PDE, their study corresponds to the model
\begin{equation}\label{june1620}
u_t=\kappa^2(t/T)u_{xx}+f(u), \  \  \  0<t<T,\ x\in\R.
\end{equation}
They proved that if $\kappa$ is increasing, and $f$ is of the Fisher-KPP type, 
the shift is algebraic and not logarithmic in time: there exists $C>0$ such that
$$
\frac{T^{1/3}}C\leq X(T)-c_{eff}T\leq CT^{1/3},~~c_{eff}=2\di\int_0^1\kappa(s)ds.
$$
In \cite{NRR1}, we proved the asymptotics
\begin{equation}\label{june1622}
X(T)=c_{eff}T-\bar\nu T^{1/3}+O({\log}T), \hbox{ with } \bar\nu=\beta\int_0^1\kappa(\tau)^{1/3}\dot\kappa(\tau)^{2/3}d\tau.
\end{equation}
Here,    $\beta < 0$ is the first zero of the Airy function. 
Maillard and Zeitouni~\cite{MZ} 
refined the asymptotics further, proving a logarithmic correction to \eqref{june1622}, and convergence of $u(T)$ to a traveling wave.

\section{Strategy of the proof of Theorem \ref{t2.10}}\label{sec:3}

\subsubsection*{Why converge to a traveling wave?}

We first provide an informal argument for the convergence of the solution of the initial
value problem to a traveling wave. Consider the Cauchy problem (\ref{e2.1}),
starting at $t=1$ for the convenience of the notation:
\begin{eqnarray}\label{sep2704}
u_t-u_{xx}=u-u^2,~~x\in\Rm,~~t>1,
%\\
%&&u(1,x)=u_1(x)=1-H(x)+v_0(x),\  \ \hbox{$v_0$ compactly supported},\nonumber
\end{eqnarray} 
 and
proceed with a standard sequence of changes of variables. 
We first go into the moving frame: 
\[
x \mapsto x - 2t + ({3}/{2}) \log t,
\]
leading to
\begin{equation}
\label{e2.0}
u_t-u_{xx}-(2-\frac3{2t})u_x=u-u^2.
\end{equation}
Next, we take out the exponential factor: set
$$
u(t,x)=e^{-x}v(t,x)
$$
so that $v$ satisfies
\begin{equation}
\label{e2.6}
v_t-v_{xx}-\frac3{2t}(v-v_x)+e^{-x}v^2=0, \quad \quad x \in \Rm, \quad t > 1.
\end{equation}
Observe that for any shift $x_\infty \in \Rm$, 
the function $V(x) = e^{x} \phi(x - x_\infty)$ is a steady solution of 
\[
V_t - V_{xx} + e^{-x} V^2 = 0.
\]
We regard (\ref{e2.6}) as a perturbation of this equation,
and expect that $v(t,x) \to e^{x} \phi(x - x_{\infty})$ as $t \to \infty$,
for some $x_\infty\in\Rm$.

\subsubsection*{The self-similar variables}

We note that for $x \to+\infty$, the term $e^{-x}v^2$ in (\ref{e2.6}) 
is negligible, while for $x\to-\infty$ the same term 
will create a large absorption and
force the solution to be close to
zero. For this reason, 
the linear Dirichlet problem
\begin{eqnarray}\label{sep1810}
&&z_t-z_{xx}-\frac3{2t}(z-z_x)=0, \quad \quad x > 0\\
&&z(t,0)=0\nonumber
\end{eqnarray}
is a reasonable proxy for (\ref{e2.6}) for $x \gg 1$, and, as shown in 
\cite{HNRR1,HNRR2}, it provides good sub- and super-solutions for $v(t,x)$. 
The main lesson of \cite{HNRR1,HNRR2} is that everything relevant to 
the solutions of~(\ref{sep1810}) 
%(and to solutions of (\ref{e2.6})) 
happens at the spatial scale~$x\sim\sqrt t$, and their asymptotics 
%of the linear problem (\ref{sep1810}) 
may be unraveled by a self-similar change of variables. Here, 
we will accept the 
full nonlinear equation \eqref{e2.6} and perform directly the 
self-similar 
change of variables
\begin{equation}
\label{e2.500}
\tau=\ln t,\ \ \ \eta=\frac{x}{\sqrt{t}}
\end{equation}
followed by a change of the unknown
$$
v(\tau,\eta)=e^{\tau/2}w(\tau,\eta).
$$
This transforms \eqref{e2.6} into
\begin{equation}
\label{e2.7}
w_\tau-\di\frac\eta2w_\eta-w_{\eta\eta}-w
+\di\frac32e^{-\tau/2}w_\eta+e^{3\tau/2-\eta{\mathrm{exp}}(\tau/2)}w^2=0,\  \   \   \eta\in\R,~~\tau>0.
\end{equation}
This transformation strengthens the reason why
the Dirichlet problem \eqref{sep1810} appears naturally: 
for
\[
\eta\ll -\tau e^{-\tau/2},
\]
the last term in the left side of (\ref{e2.7}) becomes 
exponentially large, which forces $w$ to be almost 0 in this region. 
On the other hand, for
\[
\eta\gg \tau e^{-\tau/2},
\]
this term is very small, so it 
should not play any role in the dynamics of $w$ in that region. The transition
region has width of the order $\tau e^{-\tau/2}$.

\subsubsection*{The choice of the shift}

Also, through this change of variables, we can see how a particular 
translation of the wave will be chosen. Considering (\ref{sep1810}) 
in the self-similar variables, one can show -- see \cite{HNRR1,Henderson} 
-- that, as~$\tau\to+\infty$, we have
\begin{equation}\label{sep2702}
e^{-\tau/2} z(\tau,\eta)\sim \alpha_\infty \eta e^{-\eta^2/4},~~\eta>0,
\end{equation}
with some $\alpha_\infty>0$.
Therefore, taking (\ref{sep1810}) as an approximation to (\ref{e2.6}), we should expect that
\begin{equation}\label{nov2302}
u(t,x)=e^{-x}v(t,x)\sim e^{-x}z(t,x)\sim e^{-x}e^{\tau/2}
\alpha_\infty \eta e^{-\eta^2/4}=\alpha_\infty xe^{-x}e^{-x^2/(4t)},
\end{equation}
at least for $x$ of the 
order $O(\sqrt{t})$.  This determines the unique translation: 
if we accept that $u$ converges to a translate~$x_\infty$ of $\phi_*$, then 
for large $x$ (in the moving frame) we have
\begin{equation}\label{sep1812}
u(t,x)\sim\phi_*(x - x_\infty)\sim xe^{-x + x_\infty}.
\end{equation}
Comparing this with (\ref{nov2302}), we infer that
$$
x_\infty= \ln\alpha_\infty.
$$

The difficulty with this argument, apart from the justification
of the approximation
\[
u(t,x)\sim e^{-x}z(t,x),
\]
is that each of the asymptotics (\ref{nov2302}) and (\ref{sep1812})
uses different ranges of $x$: 
(\ref{nov2302}) comes from the self-similar variables in the region
$x\sim O(\sqrt t)$, while (\ref{sep1812}) assumes $x$ to be large but finite. 
However, the self-similar analysis does not tell us at 
this stage what happens on the scale~$x\sim O(1)$. 
Indeed, it is clear from (\ref{e2.7}) that the error in the 
approximation (\ref{sep2702}) is at least of the
order~$O(e^{-\tau/2})$ -- note that the right side in (\ref{sep2702})
is a solution of (\ref{e2.7})
without the last two terms in the left side. 
On the other hand, the scale~$x\sim O(1)$ corresponds to $\eta\sim e^{-\tau/2}$.
Thus, the leading order term and the error in~(\ref{sep2702}) are of the same size
for~$x\sim O(1)$, which means that
we can not extract information directly from~(\ref{sep2702}) on that scale.

To overcome this issue, we proceed in two steps: 
first we use the self-similar variables to 
prove stabilization (that is, (\ref{nov2302}) holds) at the spatial scales $x\sim O(t^\gamma)$ with 
a small $\gamma>0$, 
and not just at the diffusive scale $O(\sqrt{t})$.  This boils down 
to showing that 
\[
w(\tau,\eta) \sim  \alpha_\infty \eta e^{-\eta^2/4}
\]
for the solution to (\ref{e2.7}), even for $\eta \sim e^{-(1/2 - \gamma) \tau}$.  
Next, we show that this stabilization 
is sufficient to ensure the stabilization on the scale $x\sim O(1)$ and convergence to a unique wave. 
This is the core of the argument: 
everything happening at $x\sim O(1)$ should be governed by the tail of the solution -- the fronts are pulled.

We conclude this section with some remarks about the generality of the argument.  Although we assume, for simplicity, that the reaction term in (\ref{e2.1}) is quadratic, our proof also works for a more general reaction term.  Specifically, the function $u - u^2$ in (\ref{e2.1}) may be replaced by a $C^2$ function $f:[0,1] \to \Rm$ satisfying $f(0) = 0 = f(1)$, $f'(0) > 0$, $f'(1) < 0$, and $f'(s) \leq f'(0)$ for all $s \in [0,1]$. In particular, these assumptions imply that there is $C > 0$ such that $0 \leq f'(0) s - f(s) \leq  C s^2$ for all $s \in [0,1]$.  Without loss of generality, we may suppose that $f'(0) = 1$. Then, if $g(u) = u - f(u)$, the equation (\ref{e2.6}) for $v$ becomes
\[
v_t-v_{xx}-\frac3{2t}(v-v_x)+e^{x}g(e^{-x} v)=0, \quad \quad x \in \Rm, \quad t > 1,
\]
and the equation (\ref{e2.7}) for $w$ becomes
\[
w_\tau-\di\frac\eta2w_\eta-w_{\eta\eta}-w +\di\frac32e^{-\tau/2}w_\eta+e^{\tau/2+\eta{\mathrm{exp}}(\tau/2)}g( e^{\tau/2-\eta{\mathrm{exp}}(\tau/2)} w) =0,\  \   \   \eta\in\R,~~\tau>0.
\]
where $0 \leq g(s) \leq C s^2$ and $g'(s) \geq 0$.  Then all of the arguments below (and in \cite{HNRR1}) work in this more general setting.  Finally, the arguments also apply to fronts arising from compactly supported initial data $u_0 \geq 0$ (not just perturbations of the step-function). In that case, one obtains two fronts propagating in opposite directions.  Combined with \cite{HNRR1}, our arguments here imply that Theorem \ref{t2.10} holds for both fronts. That is, the fronts moving to $\pm \infty$ are at positions $\sigma_\infty^\pm(t)$ with
\[
\sigma^\pm_\infty(t) = \pm 2 \mp \frac{3}{2} \log t + x_\infty^{\pm} + o(1)
\]
where the shifts $x_\infty^+$ and $x_\infty^-$ may differ and depend on the initial data.

\section{Convergence to a single wave as a consequence of the diffusive scale convergence}\label{sec:4}

The proof of Theorem~\ref{t2.10} relies on the following two lemmas.
The first is a consequence of~\cite{HNRR1}.
\begin{lem}\label{lem:mar2302}
The solution of (\ref{e2.0}) with $u(1,x) = u_0(x)$ satisfies
\begin{equation}\label{mar2306}
\lim_{x\to-\infty}u(t,x)=1,~~\lim_{x\to+\infty}u(t,x)=0,
\end{equation}
both uniformly in $t>1$.
\end{lem}
The main new step is to establish 
the following.
\begin{lem}\label{sep29-lem02}
There exists a constant $\alpha_\infty>0$ with the following property.
For any $\gamma>0$ and all $\eps>0$ we can find $T_\eps$
so that for all $t>T_\eps$ we have 
\begin{equation}\label{sep2902}
|u(t,x_\gamma)-\alpha_\infty x_\gamma e^{-x_\gamma}e^{-x_\gamma^2/(4t)}|\le
\eps x_\gamma e^{-x_\gamma}e^{-x_\gamma^2/(6t)},
\end{equation}
with $x_\gamma=t^\gamma$.
\end{lem}
We postpone the proof of this lemma for the moment, and show how it is used.
A consequence of Lemma~\ref{sep29-lem02} is that the problem for the moment
is to understand, for a given $\alpha>0$,
the behavior of the solutions of
\begin{eqnarray}\label{e4.4}
&&\pdr{u_\alpha}{t}-\frac{\partial^2u_\alpha}{\partial x^2}
-\di(2-\di\frac3{2t})\pdr{u_\alpha}{x}-u_\alpha+u_\alpha^2=0,
\  \  x\leq x_{\gamma}(t)\\
&&u_\alpha(t,t^\gamma)=\alpha t^\gamma e^{-t^\gamma-t^{2\gamma-1}/4},\nonumber
\end{eqnarray}
for $t>T_\eps$, with the initial condition $u_\alpha(T_\eps,x)=u(T_\eps,x)$.
%Then we will consider the two solutions
%``inner" problems
%\begin{eqnarray}\label{sep2534}
%&&\tilde u_t-\tilde u_{xx}-\di(2-\di\frac3{2t})u_x-u+u^2=0\  \  x\leq x_\gamma (t)\\
%&&u(t,x_\gamma)=(\alpha_\infty\pm\e)
%x_\gamma e^{-x_\gamma -x_\gamma^{2}/4t},\nonumber
%\end{eqnarray}
%motivated by the error term in (\ref{sep2902}),and 
%Taking $\alpha=\alpha_\infty\pm\eps$, 
In particular, we will show that $u_{\alpha_\infty\pm\eps}(t,x)$
converge, as $t\to+\infty$, to a pair of steady solutions, separated only 
by an order~$O(\e)$-translation.
Note that the function $v(t,x)=e^{x}u_\alpha(t,x)$   solves
\begin{eqnarray}
\label{e4.1}
&&v_t-v_{xx}+\di\frac3{2t}(v_x-v)+e^{-x}v^2=0,\  \  x\leq t^{\gamma}\\
&&v(t,t^\gamma)=\alpha t^\gamma e^{-t^{2\gamma-1}/4}.\nonumber
\end{eqnarray}
Since we anticipate that the tail is going to dictate the behavior of $u_\alpha$, 
we choose the translate of the wave 
that matches exactly the behavior of $u_\alpha(t,x)$ 
at the boundary $x=t^\gamma$: set 
\begin{equation}
\label{e4.3}
\psi(t,x)=e^{x}\phi_*(x+\zeta(t)).
\end{equation}
Recall that $\phi_*(x)$ is the traveling wave profile.
We look for a function $\zeta(t)$ in (\ref{e4.3}) such that 
\begin{equation}
\label{e4.2}
\psi(t,t^\gamma)=v(t,t^\gamma).
\end{equation}
In view of the expansion \eqref{e2.20}, we should have, with some $\omega_0>0$:
$$
e^{-\zeta(t)}(t^\gamma+\zeta(t)+k)+O(e^{-\omega_0 t^\gamma})=
\alpha t^\gamma e^{-1/(4t^{1-2\gamma})},
$$
which implies, for $\gamma \in (0,1/2)$,
\[
\zeta(t)=-\ln \alpha-(\ln\alpha-k)t^{-\gamma}+O(t^{-2\gamma}), 
\]
and thus
\[
|\dot \zeta(t)|\le\di\farc{C}{t^{1+\gamma}}.
\]
The equation for the function $\psi$ is 
$$
\psi_t-\psi_{xx}+\di\frac3{2t}(\psi_x-\psi)+e^{-x}\psi^2=
- \dot \zeta \psi+ \dot\zeta \psi_{x}+\di\frac3{2t}(\psi_x-\psi)=
O(\frac{x}t)=O(t^{-1+\gamma}),~~|x|<t^\gamma.
$$
In addition, the left side above is exponentially small for $x<-t^{\gamma}$
because of the exponential factor in (\ref{e4.3}).
Hence, the difference $
s(t,x)=v(t,x)-\psi(t,x)$
satisfies
\begin{eqnarray}
\label{e4.5}
&&s_t-s_{xx}+\di\frac3{2t}(s_x-s)+e^{-x}(v+\psi)s=O(t^{-1+\gamma}),\  \  \vert x\vert\leq t^{\gamma}\\
&&s(t,-t^{\gamma})=O(e^{-t^\gamma}),\   \   s(t,t^\gamma)=0.\nonumber
\end{eqnarray}
\begin{prop}
\label{p4.1}
For $\gamma \in (0,1/3)$, we have 
\begin{equation}
\label{mar2210}
\di\lim_{t\to+\infty}\sup_{\vert x\vert \leq t^\gamma}\vert s(t,x)\vert=0.
\end{equation}
\end{prop}
\noindent {\bf Proof.} 
%The Kato inequality transforms equation \eqref{e4.5} into the linear inequation
%\begin{equation}
%\label{e4.6}
%\begin{array}{rll}
%\vert w\vert_t-\vert w\vert_{xx}+\di\frac3{2t}(\vert w\vert_x-\vert w\vert)+e^{-x}(v+\psi)w=&O(t^{-1+\gamma}),\  \  \vert x\vert\leq t^{\gamma}\\
%\vert w(t,-t^{\gamma})\vert=O(e^{-t^\gamma}),\   \   \vert w\vert(t,t^\gamma)=&0.
%\end{array}
%\end{equation}
The issue is whether the Dirichlet boundary conditions would be stronger
than the force in the right side of (\ref{e4.5}).  Since
the principal Dirichlet eigenvalue for the Laplacian 
in $(-t^{\gamma},t^{\gamma})$ is~$\frac{\pi^2}{4t^{2\gamma}}$, 
investigating \eqref{e4.5} is, heuristically, equivalent to solving the ODE
\begin{equation}\label{mar2208}
f'(t)+(1-2\gamma){t^{-2\gamma}} f=\frac1{t^{1-\gamma}}.
\end{equation}
The coefficient $(1-2\gamma)$ is chosen simply for convenience and
can be replaced by another constant. The solution of (\ref{mar2208}) is
%\[
%\frac{d}{dt}(e^{t^{-2\gamma+1} }f)=
%t^{3\gamma-1}e^{t^{-2\gamma+1} },
%\]
%so that
\[
f(t)=f(1)e^{(-t^{-2\gamma+1}+1) }+\int_1^t
s^{\gamma-1}e^{(-t^{-2\gamma+1}+s^{-2\gamma+1}) }ds.
\]
Note that $f(t)$ tends to 0 as $t\to+\infty$
a little faster than $t^{3\gamma-1}$ as soon as $\gamma<1/3$, so the analog
of~(\ref{mar2210}) holds for the solutions of (\ref{mar2208}).
With this idea in mind, we are going to look for a super-solution
to \eqref{e4.5}, in the form
\begin{equation}
\label{e4.7}
\overline s(t,x)=t^{-\lambda}\cos \left(\frac{x}{t^{\gamma+\e}} \right),
\end{equation}
where $\lambda$, $\gamma$ and $\e$ will be chosen to be small enough. 
We now set $T_\eps=1$ for convenience.
We have, for~$\vert x\vert\leq t^{\gamma}$:
\begin{eqnarray}
\label{e4.8}
&&\overline s(t,x)\sim t^{-\lambda},\  
-\overline s_{xx}= t^{-(2\gamma+2\e)}\overline s(t,x),\\
&&\overline s_t=-\farc{\lambda}{t}
\bar s+g(t,x),~~|g(t,x)|\le\farc{C|x|}{t^{\lambda+\gamma+\eps+1}}\le 
\farc{C}{t^{1+\eps}}\overline s(t,x),\nonumber
\end{eqnarray}
and 
\begin{equation}
\label{e4.9}
\left \lvert \di\frac3{2t} (\overline s_x-\overline s)(t,x) \right \rvert \leq Ct^{-1}\overline s(t,x).
\end{equation}
Gathering \eqref{e4.8} and \eqref{e4.9} we infer the 
existence of $q>0$ such that, for $t$ large enough:
$$
\biggl(\partial_t-\partial_{xx}+\di\frac3{2t}(\partial_x-1)\biggl)\overline s(t,x)
\geq qt^{-(2\gamma+2\e)}\overline s(t,x)\geq 
\frac{q}2t^{-(2\gamma+2\e+\lambda)}\geq O(\frac1{t^{1-\gamma}}),
$$
as soon as $\e$ and $\lambda$ are small enough, since $\gamma \in (0,1/3)$. Because 
the right side of \eqref{e4.5} does not depend on $\overline s$, 
the inequality extends to all $t\ge 1$ 
by replacing $\overline s$ by $A\overline s$, with $A$ large enough,
and (\ref{mar2210}) follows.

Let us note that the term $e^{-x}(v + \psi)$ in (\ref{e4.5}), which results from the quadratic structure of the nonlinearity, is positive.  For a more general nonlinearity $f(u)$ replacing $u - u^2$, the monotonicity of $g(u) = uf'(0) - f(u)$ may be used in an analogous way.  ~$\Box$

\subsubsection*{Proof of Theorem \ref{e2.1}}

We are now ready to prove the theorem. Fix $\gamma \in (0,1/3)$, as required by Proposition \ref{p4.1}. Given $\eps>0$, take
$T_\eps$ as in Lemma~\ref{sep29-lem02}. Let $u_\alpha(t,x)$ be the solution of~\eqref{e4.4} for $t>T_\eps$,
and the initial condition $u_\alpha(T_\eps,x)=u(T_\eps,x)$. Here,
$u(t,x)$ is the solution of the original problem \eqref{e2.0}. Taking $T_\epsilon$ larger, if necessary, we may assume that $e^{-x_\gamma^2/(4t)} \geq 1/2$ for $t \geq T_\epsilon$. It follows from Lemma~\ref{sep29-lem02} that
for any $t\geq T_\e$, we have 
\[
u_{\alpha_\infty-2\e}(t,x)\leq u(t,x)\leq u_{\alpha_\infty+2\e}(t,x),
\]
for all $x\le t^\gamma$. 
From Proposition \ref{p4.1}, we have 
\begin{equation}\label{mar2302}
e^x\big[u_{\alpha_\infty\pm 2\e}(t,x)-\phi_*(x+\zeta_\pm(t))\big]=o(1),
\hbox{ as $t\to+\infty$},
\end{equation}
uniformly in $x\in(-t^\gamma,t^\gamma)$, with
$$
\zeta_\pm(t)=-(1-t^{-\gamma})\ln(\alpha_\infty\pm 2\e)+O(t^{-2\gamma}).
$$
Because $\e>0$ is arbitrary, we have  
$$
\lim_{t\to+\infty}\big(u(t,x)-\phi_*(x+x_\infty)\big)=0,
$$
with $x_\infty=-\ln\alpha_\infty$, uniformly on compact sets.
Together with Lemma~\ref{lem:mar2302},
this concludes the proof of Theorem \ref{t2.10}.~$\Box$

\section{The diffusive scale $x\sim O(\sqrt{t})$ 
and the proof of Lemma~\ref{sep29-lem02}}\label{sec:5}

Our analysis starts with (\ref{e2.7}),
which we write as
\begin{equation}
\label{sep2102}
w_\tau+Lw
+\di\frac32e^{-\tau/2}w_\eta+e^{3\tau/2-\eta{\mathrm{exp}}(\tau/2)}w^2=0,\  \   \   \eta\in\R,~~\tau>0.
\end{equation}
Here, the operator $L$ is defined as
\begin{equation}\label{sep2104}
Lv=-v_{\eta\eta}-\di\frac{\eta}{2}v_\eta-v.
\end{equation}
Its principal eigenfunction 
on the half-line $\eta>0$ with the Dirichlet boundary condition at~$\eta=0$~is 
\[
\phi_0(\eta)=\farc{\eta}{2} e^{-\eta^2/4},
\]
as $L\phi_0=0$. The operator $L$ has a discrete spectrum in $L^2(\R_+)$, weighted by 
$e^{-\eta^2/8}$, its non-zero eigenvalues are $\lambda_k = k \ge 1$, and the 
corresponding eigenfunctions are related via
\[
\phi_{k+1}=\phi_k''.
\]
The principal eigenfunction of the adjoint operator 
\[
L^*\psi=-\psi_{\eta\eta}+\farc 12\partial_\eta(\eta \psi)-\psi
\]
is~$\psi_0(\eta)=\eta$. Thus, the solution of the unperturbed version of (\ref{sep2102})
on a half-line
\begin{equation}
\label{sep2106}
p_\tau+Lp=0,~~~\eta>0,~p(\tau,0)=0,
\end{equation}
satisfies
\begin{equation}\label{sep2108}
p(\tau,\eta)=\eta\ \frac{e^{-\eta^2/4}}{2\sqrt{\pi}}
\int_0^{+\infty}\xi v_0(\xi)d\xi+O(e^{-\tau})e^{-\eta^2/6},
\hbox{ as $\tau\to +\infty$},
\end{equation}
and our task is to generalize  this asymptotics to the full problem 
(\ref{sep2102}) on the whole line.  The weight~$e^{-\eta^2/6}$ in (\ref{sep2108})
is, of course, by no means optimal. We will prove the following:

\begin{lem}
\label{l3.1bis}
Let $w(\tau,\eta)$ be the solution of \eqref{e2.7} on $\Rm$, with 
the initial condition $w(0,\eta)=w_0(\eta)$ 
such that~$w_0(\eta)=0$ for all $\eta>M$, with some $M>0$, and $w_0(\eta)=O(e^{\eta})$
for $\eta<0$. 
There exists~$\alpha_\infty>0$ and a function $h(\tau)$ such that 
$
\di\lim_{\tau\to+\infty}h(\tau)=0,$
and such that we have, for any $\gamma'\in(0,1/2)$:
\begin{equation}\label{mar2312}
w(\tau,\eta)=(\alpha_\infty+h(\tau))\eta_+ 
e^{-\eta^2/4}+ R(\tau,\eta)e^{-\eta^2/6}, \quad \quad \eta \in \Rm,
\end{equation}
with 
$$|R(\tau,\eta)|\le C_{\gamma'}e^{-(1/2-\gamma')\tau},
$$
and where $\eta_+ = \max(0,\eta)$.
\end{lem}
Once again, the weight $e^{-\eta^2/6}$ is not optimal.
Lemma~\ref{sep29-lem02} is an immediate consequence of this result. 
Indeed,  
\[
u(t,x)=e^{-x}\sqrt{t}w(\log t,\frac{x}{\sqrt{t}}),
\]
hence Lemma~\ref{l3.1bis} implies, with $x_\gamma=t^\gamma$,
\begin{eqnarray}\label{mar2314}
&&e^{x_\gamma}u(t,x_\gamma)-\alpha_\infty x_\gamma e^{-x_\gamma^2/(4t)}
= \sqrt{t}w\Big(\log t,\farc{x_\gamma}{\sqrt{t}}\Big)
-\alpha_\infty x_\gamma e^{-x_\gamma^2/(4t)}\\
&&~~~~~~~~~~~~~~=
h(\log t)x_\gamma e^{-x_\gamma^2/(4t)}+\sqrt{t}
R\Big(\log t,\farc{x_\gamma}{\sqrt{t}}\Big)e^{-x_\gamma^2/(6t)}.\nonumber
\end{eqnarray}
We now take $T_\eps$ so that $|h(\log t)|<\eps/3$
for all $t>T_\eps$. For the second term in the right side of (\ref{mar2314})
we write
\begin{eqnarray}\label{mar2316}
\big|R\big(\log t,\farc{x_\gamma}{\sqrt{t}}\big)\big|\sqrt{t}
e^{-x_\gamma^2/(6t)}\le C t^{\gamma'}e^{-x_\gamma^2/(6t)}\le 
\eps x_\gamma e^{-x_\gamma^2/(6t)}
\end{eqnarray}
for $t>T_\eps$ sufficiently large, as soon as $\gamma'<\gamma$.
This proves (\ref{sep2902}). Thus, the proof of Lemma~\ref{sep29-lem02} 
reduces to proving Lemma~\ref{l3.1bis}. We will prove the latter by a construction of an upper and lower barrier for
$w$ with the correct behaviors.

\subsubsection*{The approximate Dirichlet boundary condition}

Let us come back to
 why the solution of (\ref{sep2102}) must approximately satisfy the Dirichlet boundary condition
at $\eta=0$. 
Recall that $w$ is related to the solution of the original KPP problem via
\[
w(\tau,\eta)=u(e^\tau,\eta e^{\tau/2})e^{-\tau/2+\eta e^{\tau/2}}.
\]
%we have
%\[
%w(\tau,-e^{-\tau/3})=u(e^\tau-1,-e^{\tau/6})e^{-\tau/2- e^{\tau/6}}.
%\]
The trivial a priori bound $0<u(t,x)< 1$ implies that  we have
\begin{equation}\label{sep2120}
0<w(\tau,\eta)<e^{-\tau/2+\eta e^{\tau/2}},~~\eta<0,
\end{equation}
and, in particular, we have
\begin{equation}\label{sep2110}
0<w(\tau,-e^{-(1/2-\gamma)\tau})\le e^{-e^{\gamma\tau}}.
\end{equation}
We also have
\[
w_\tau(\tau,\eta)= u_t(e^\tau,\eta e^{\tau/2})e^{\tau/2+\eta e^{\tau/2}}+\frac\eta2 u_x(e^\tau,\eta e^{\tau/2})e^{\eta e^{\tau/2}}
+(\frac\eta2 e^{\tau/2}-\farc{1}{2})u(e^\tau,\eta e^{\tau/2})e^{-\tau/2+\eta e^{\tau/2}},
\]
so that
\begin{eqnarray}\label{sep2114}
&&w_\tau(\tau,-e^{-(1/2-\gamma)\tau})=u_t(e^\tau,-e^{\gamma \tau})
e^{\tau/2- e^{\gamma\tau}}-
\farc 12
e^{-(1/2-\gamma)\tau}u_x(e^\tau,-e^{\gamma \tau})e^{-e^{\gamma\tau}}
\nonumber\\
&&~~~~~~~~~~~~~
~~~~~~~~~~~~-\frac12( e^{\gamma\tau}+1)u(e^\tau,- e^{\gamma\tau})e^{-\tau/2- e^{\gamma\tau}}\nonumber\\
&&~~~~~~~~~~~~~
~~~~~~~~~~~~=O(e^{-\gamma e^{\gamma\tau}}),
\end{eqnarray}
for $\gamma>0$ sufficiently small. 
Thus, the solution of (\ref{sep2102}) satisfies
\begin{equation}\label{sep2116}
\begin{array}{rll}
0<w(\tau,-e^{-(1/2-\gamma)\tau})\le &e^{-e^{\gamma\tau}},\\
|w_\tau(\tau,-e^{-(1/2-\gamma)\tau})|\le &Ce^{-\gamma e^{\gamma\tau}},
\end{array}
\end{equation}
which  we will use as an approximate Dirichlet boundary condition at $\eta=0$. 
 
\subsubsection*{An upper barrier} 
 
Consider the solution of 
\begin{eqnarray}\label{sep2122}
&&\overline w_\tau+L\overline w+\di\frac32
e^{-\tau/2}\overline w_\eta=0,~~\tau>0,~~\eta>-e^{-(1/2-\gamma)\tau},\\
&&  \overline w(\tau,-e^{-(1/2-\gamma)\tau})= e^{-e^{\gamma\tau}},\nonumber
\end{eqnarray}
%with $C$ large enough so that 
%\[
%\overline w(\tau, -e^{-\tau/3})\geq w(\tau, -e^{-\tau/3}),
%\]
%and 
with a compactly supported  initial condition 
$\bar w_0(\eta)=\bar w(0,\eta)$ chosen so that 
$\bar w_0(\eta)\ge u(1,\eta)e^{\eta}.$ Here,~$\gamma\in(0,1/2)$ should be thought of as a small parameter.

It follows from (\ref{sep2116}) that $\overline w(\tau,\eta)$
is an upper barrier for $w(\tau,\eta)$.
That is, we have
\[
w(\tau,\eta)\le \bar w(\tau,\eta),
\hbox{ for all $\tau>0$ and $\eta>-e^{-(1/2-\gamma)\tau}$}.
\]
It is convenient to make a change of variables 
\begin{equation}
\bar w(\tau,\eta)=\bar p(\tau,\eta+e^{-(1/2-\gamma)\tau})+e^{-e^{\gamma\tau}}g(\eta+e^{-(1/2-\gamma)\tau}),
\end{equation}
where $g(\eta)$ is a smooth monotonic function such that  $g(\eta)=1$ for $0\le\eta<1$ and $g(\eta)=0$ for $\eta>2$.
The function $\bar p$ satisfies
\begin{equation}\label{sep2124}
\bar p_\tau+L\bar p+(\gamma e^{-(1/2-\gamma)\tau}+\farc32 
e^{-\tau/2})\bar p_\eta=G(\tau,\eta)e^{- e^{\gamma\tau}},~~\eta>0,~~\bar p(\tau,0)=0,
\end{equation}
for $\tau > 0$, with a smooth function $G(\tau,\eta)$ supported in $0\le\eta\le 2$, and the initial condition
\[
\bar p_0(\eta)=\bar w_0(\eta-1)-e^{-1}g(\eta),
\]
which also is compactly supported. 

We will allow (\ref{sep2124}) to run for a large time $T$, after which time
we can treat the right side and the last term in the left side of (\ref{sep2124})
as a small perturbation. A variant of Lemma 2.2 from~\cite{HNRR1} implies that $\bar p(T,\eta) e^{\eta^2/6} \in L^2(\Rm_+)$ for all $T > 0$, as well as the following estimate:

\begin{lem}\label{l3.2bis}
Consider $\omega\in(0,1/2)$ and $G(\tau,\eta)$ smooth, bounded, and compactly supported in $\R_+$. Let $p(\tau,\eta)$ solve
\begin{equation}
\label{e3.2}
\vert p_\tau+Lp\vert\leq \epsilon 
e^{-\omega \tau}(\vert p_\eta\vert+\vert p\vert+G(\tau,\eta)),\ \ 
\tau>0,\,\eta>0,\ \ \ \ \ \ \ \ p(\tau,0)=0.
\end{equation}
with the initial condition $p_0(\eta)$ such that
 $p_0(\eta) e^{\eta^2/6} \in L^2(\Rm_+)$. There exists $\epsilon_0>0$ and $C>0$ (depending on $p_0$) such that, for all~$0<\epsilon<\epsilon_0$, we have
\begin{equation}\label{sep2131}
p(\tau,\eta)=\eta\biggl(\frac{e^{-\eta^2/4}}{2\sqrt{\pi}}
\Big(\int_0^{+\infty}\xi p_0(\xi)d\xi+\epsilon R_1(\tau,\eta)\Big)+ 
\e e^{-\omega\tau}R_2(\tau,\eta)e^{-\eta^2/6}+e^{-\tau}R_3(\tau,\eta)e^{-\eta^2/6}\biggl),
\end{equation}
where $\|R_{1,2,3}(\tau,\cdot)\|_{C^3}\leq C $ for all $\tau>0$.
\end{lem}

For any $\epsilon > 0$, we may choose $T$ sufficiently large, and $\omega \in (0, 1/2-\gamma)$ so that
\begin{equation}
\label{peqnabs}
\vert \bar p_\tau+L \bar p\vert\leq \epsilon 
e^{-\omega (\tau - T)}(\vert \bar p_\eta\vert+|G(\tau,\eta)|),\ \ 
\tau>T,\,\eta>0,\ \ \ \ \ \ \ \ p(\tau,0)=0.
\end{equation}
This follows from (\ref{sep2124}). Then, applying Lemma \ref{l3.2bis} for $\tau > T$, we have 
\begin{equation}\label{sep2502}
\bar p(\tau,\eta)=\eta\biggl(\frac{e^{-\eta^2/4}}{2\sqrt{\pi}}
\Big(\int_0^{+\infty}\!\!\xi \bar p(T,\xi)d\xi+\epsilon R_1(\tau,\eta)\Big)+ 
\e e^{-\omega(\tau-T)}R_2(\tau,\eta)e^{-\eta^2/6}
+e^{-(\tau-T)}R_3(\tau,\eta)e^{-\eta^2/6}\biggl).
\end{equation}

We claim that with a suitable choice of $\bar w_0$, the integral term in (\ref{sep2502}) is bounded from below: 
\begin{equation}\label{sep2130}
\int_0^\infty \eta \bar p(\tau,\eta)d\eta\ge  1, \hbox{ for all $\tau>0$.}
\end{equation}
Indeed, multiplying (\ref{sep2124}) by $\eta$ and integrating gives
\begin{eqnarray}\label{sep2126}
\frac{d}{d\tau}\int_0^\infty \eta\bar p(\tau,\eta)d\eta=(\gamma e^{-(1/2-\gamma)\tau}+\farc 32 e^{-\tau/2})\int_0^\infty \bar p(\tau,\eta)d\eta+
e^{-e^{\gamma\tau}}\int G(\tau,\eta)\eta d\eta.
\end{eqnarray}
The function $G(\tau,\eta)$ need not have a sign, hence a priori we do not know that $\bar p(\tau,\eta)$
is positive everywhere. However,
it follows from (\ref{sep2124}) that the negative part of $\bar p$ is bounded
as
\[
\di\int_0^\infty \bar p(\tau,\eta)d\eta \geq - C_0,
\]
for all $\tau>0$, with the constant  $C_0$ which does not depend on $\bar w_0(\eta)$ on the interval $[2,\infty)$. 
Thus, we deduce from (\ref{sep2126}) that for all $\tau>0$ we have
\begin{equation}\label{sep2128}
\int_0^\infty \eta \bar p(\tau,\eta)d\eta\ge  \int_0^\infty \eta\bar w_0(\eta)d\eta-C_0',
\end{equation}
with, once again, $C_0'$ independent of $\bar w_0$. Therefore, after possibly increasing $\bar w_0$ we may ensure that~(\ref{sep2130}) holds.

It follows from (\ref{sep2130}) and (\ref{sep2502}) that there exists a 
sequence $\tau_n\to+\infty$, $C>0$ and a function~$\overline W_\infty(\eta)$
such that 
\begin{equation}\label{sep2504}
C^{-1} \eta e^{-\eta^2/4}\le \overline W_\infty(\eta)\le C\eta  e^{-\eta^2/4},
\end{equation}
and
\be
\lim_{n\to+\infty} e^{\eta^2/8}\vert\bar p(\tau_n,\eta)-\overline W_\infty(\eta)\vert=0,
\ee
uniformly in $\eta$ on the half-line $\eta\ge 0$. The same bound for the function
$\bar w(\tau,\eta)$ itself follows: 
\begin{equation}\label{sep2506}
\lim_{n\to+\infty} e^{\eta^2/8}\vert\bar w(\tau_n,\eta)-\overline W_\infty(\eta)\vert=0,
\end{equation}
also uniformly in $\eta$ on the half-line $\eta\ge 0$.

\subsubsection*{A lower barrier}

A lower barrier for $w(\tau,\eta)$ is devised as follows. First, note that the 
upper barrier for $w(\tau,\eta)$ we have constructed above implies that
\[
e^{3\tau/2-\eta{\mathrm{exp}}(\tau/2)}w(\tau,\eta)\leq 
C_\gamma e^{- {\mathrm{exp}(\gamma\tau/2)}},
\]
as soon as 
\[
\eta\geq  e^{-(1/2-\gamma)\tau},
\]
with $\gamma\in(0,1/2)$,
and $C_\gamma >0$ is chosen sufficiently large. 
Thus, a lower barrier $\underline w(\tau,\eta)$ can be defined as the solution of
\begin{equation}\label{sep2508}
\underline w_\tau+L\underline w+\di\frac32
e^{-\tau/2}\underline w_\eta+C_\gamma e^{-{\mathrm{exp}(\gamma\tau/2)}}\underline w=0,\ \ \ \underline 
w(\tau,e^{-(1/2-\gamma)\tau})=0,~~~\eta> e^{-(1/2-\gamma)\tau},
\end{equation}
and with an initial condition $\underline w_0(\eta)\le w_0(\eta)$.   This time it is convenient to make the change of variables
$$
\underline w(\tau,\eta)={\underline z}(\tau,\eta-e^{-(1/2-\gamma)\tau})
$$ 
so that
\begin{equation}
\label{april2116}
\underline z_\tau+L\underline z+(-\gamma e^{-(1/2-\gamma)\tau}+\farc32 
e^{-\tau/2})\underline z_\eta + C_\gamma e^{-{\mathrm{exp}(\gamma\tau/2)}} \underline z =0,~~\eta>0,~~\underline z(\tau,0)=0,
\end{equation}
We could now try to use an abstract stable manifold theorem to prove that 
\begin{equation}\label{may202}
\underline I(\tau):=\di\int_{0}^\infty\eta\underline z(\tau,\eta)d\eta\ge c_0 > 0,~~\hbox{ for all $\tau>0$}.
\end{equation}
That is, $\underline I(\tau)$ 
remains uniformly bounded away from 0. However, to keep this paper self-contained, we give a direct proof of (\ref{may202}). We look for a sub-solution to (\ref{april2116}) in the form
\begin{equation}
\label{e5.3000}
\underline p(\tau,\eta)= \left( \zeta(\tau)\phi_0(\eta)-q(\tau)\eta e^{-\eta^2/8} \right) e^{-F(\tau)},
\end{equation}
where 
\[
F(\tau) = \int_0^\tau C_\gamma e^{- \exp(\gamma s/2)} \,ds,
\]
and with the functions $\zeta(\tau)$ and $q(\tau)$ satisfying
\begin{equation}\label{may204}
\zeta(\tau)\geq\zeta_0>0,\quad \dot\zeta(\tau)<0,\quad q(\tau)>0,\quad q(\tau)=O(e^{-\tau/4}).
\end{equation}
In other words,  we wish to devise $\underline p(\tau,\eta)$ as in (\ref{e5.3000})-(\ref{may204})
such that 
\begin{equation}\label{may208}
\underline p(0,\eta)\le \underline z(0,\eta)=w_0(\eta+1),
\end{equation}
and
\begin{equation}\label{may206}
{\cal L}(\tau)\underline p\leq0,
\end{equation}
with
$$
{\cal L}(\tau)\underline p= \underline p_\tau+L\underline p+(-\gamma e^{-(1/2-\gamma)\tau}+\farc32 
e^{-\tau/2})\underline p_\eta.
$$
Notice that the choice of $F(\tau)$ in (\ref{e5.3000}) has eliminated a low order term involving $C_\gamma e^{- \exp(\gamma \tau/2)}$. For convenience, let us define 
\[
h(\tau) = -\gamma e^{-(1/2-\gamma)\tau}+\frac{3}{2}
e^{-\tau/2},
\]
which appears in (\ref{april2116}).  Because $L\phi_0=0$ and because
$$
L(\eta e^{-\eta^2/8}) = \eta L e^{-\eta^2/8} =(\frac{\eta^2}{16}-\frac{3}{4})\eta e^{-\eta^2/8},
$$
we find that
\begin{eqnarray}\label{may210} 
{\cal L}(\tau)\underline p=\dot\zeta\phi_0+\zeta h(\tau) \phi_0' -\biggl(\dot q+(\frac{\eta^2}{16}-\frac{3}{4})q\biggl)\eta e^{-\eta^2/8} + q\frac{\eta^2}{4} e^{-\eta^2/8} h(\tau)  - q e^{-\eta^2/8} h(\tau). \nonumber
\end{eqnarray}
Let us write this as
\begin{eqnarray}\label{Lequivalent}
\eta^{-1} e^{\eta^2/8} {\cal L}(\tau)\underline p=\dot \zeta \eta^{-1} \phi_0e^{\eta^2/8} + \eta^{-1} h(\tau)\left( \zeta e^{\eta^2/8} \phi_0' + q \left(\frac{\eta^2}{4} - 1 \right)\right)  -\biggl(\dot q+(\frac{\eta^2}{16}-\frac{3}{4})q\biggl) . 
\end{eqnarray}
Our goal is to choose $\zeta(\tau)$ and $q(\tau)$ such that (\ref{may204}) holds and the right side of (\ref{Lequivalent}) is non-positive after a certain time $\tau_0$, possibly quite large. However, and this is an important point, this time $\tau_0$ will not depend on the initial condition 
$w_0(\eta)$. 

Let us restrict the small parameter $\gamma$ to the interval $(0,1/4)$.  Observe that if $\tau_0 > 0$ is sufficiently large, 
then $h(\tau) < 0$ and $|h(\tau)| \leq e^{-\tau/4}$ for all $\tau \geq \tau_0$.  As $\phi_0(\eta) = \eta e^{-\eta^2/4}$, note that 
in (\ref{Lequivalent}) both~$\phi_0'(\eta)e^{\eta^2/8}$ and~$\phi_0(\eta)e^{\eta^2/8}$  
are bounded functions. In particular, if $\tau_0$ is large enough then
\[
|\phi_0' e^{\eta^2/8} h(\tau)| \leq e^{-\tau/4}
\]
for all $\tau \geq \tau_0$, $ \eta \geq 0$.
%Hence, there is $\eta_\gamma > 0$ so that if $\tau_0$ is large enough, then 
%\[
%\left| h(\tau) \phi_0'e^{\eta^2/8}+C_\gamma e^{-{\mathrm{exp}(\gamma\tau/2)}}\phi_0 e^{\eta^2/8} \right|  \leq \frac{1}{2} h(\tau) \phi_0'(\eta) e^{\eta^2/8} \leq  e^{-\tau/4}
%\]
%holds for all $\eta \geq \eta_\gamma$ and $\tau \geq \tau_0$. 
%\begin{equation}
%\label{e5.3004}
%\farc14(\gamma e^{-(1/2-\gamma)\tau}+\farc32 
%e^{-\tau/2})\vert\phi_0'(\eta)\vert e^{\eta^2/8}+C_\gamma e^{-{\mathrm{exp}(\gamma\tau/2)}}\phi_0(\eta)e^{\eta^2/8}\leq C_\gamma e^{-\tau/4}.
%\end{equation}

Note also that for all $\eta \geq \eta_1 = \sqrt{28}$ we have
\begin{equation}
\label{e5.3001}
\frac{\eta^2}{16}- \frac{3}{4} \geq 1 \quad \quad \text{and} \quad \quad \frac{\eta^2}{4}  - 1 \geq 0.
\end{equation}
Therefore, on the interval $\eta \in [\eta_1,\infty)$ and for $\tau \geq \tau_0$, (\ref{Lequivalent}) is bounded by
\begin{eqnarray}
\eta^{-1} e^{\eta^2/8} {\cal L}(\tau)\underline p & \leq & \eta^{-1} h(\tau) \zeta e^{\eta^2/8} \phi_0'  -\left(\dot q+q\right) \leq \zeta(\tau) e^{-\tau/4} -\left(\dot q+q\right), \nonumber
\end{eqnarray}
assuming $q(\tau) > 0$ and $\dot \zeta < 0$.  Hence, if $q(\tau)$ and $\zeta(\tau)$ are chosen to satisfy 
the differential inequality
\begin{equation}
\label{e5.3006} 
\dot q+q-  e^{-\tau/4}\zeta\geq0,\quad \tau\geq\tau_0,
\end{equation}
then we will have 
\begin{equation}\label{may212}
\hbox{${\cal L}(\tau)\underline p\leq 0$ for $\tau\geq\tau_0$ and $\eta\geq\eta_1$,}
\end{equation}
provided that $\dot\zeta\leq0$, as presumed in (\ref{may204}). Still 
assuming $\dot\zeta\leq 0$ on $(\tau_0,+\infty)$, a sufficient condition for (\ref{e5.3006}) to be satisfied is:
$$
\dot q+q\geq e^{-\tau/4}\zeta(\tau_0),\quad \tau\geq\tau_0.
$$
Hence, we choose
\begin{equation}
\label{e5.3007} 
q(\tau)= e^{-(\tau-\tau_0)}+\frac{4 }3e^{-\tau /4}\zeta(\tau_0).
\end{equation}
Note that $q(\tau)$ satisfies the assumptions on $q$ in (\ref{may204}).
%So, enlarging $C_\gamma$ we choose:
%\begin{equation}
%\label{e5.3007} 
%q(\tau)= q(\tau_0)e^{-(\tau-\tau_0)}+e^{-\tau_0/4}C_\gamma e^{-(\tau-\tau_0)/4}\zeta(\tau_0).
%\end{equation}

Let us now deal with the range $\eta\in[0,\eta_1]$. 
The function $\eta^{-1}\phi_0(\eta)$ is bounded on $\Rm$ and it is bounded away from $0$ on $[0,\eta_1]$. Define
\[
\epsilon_1 = \min_{\eta \in [0,\eta_\gamma]} \eta^{-1} \phi_0(\eta) e^{\eta^2/8} > 0.
\]
As $h(\tau)<0$ for $\tau\ge\tau_0$,
on the interval $[0,\eta_1]$, we can bound (\ref{Lequivalent}) by
\begin{eqnarray}\label{Lequivalent2}
\eta^{-1} e^{\eta^2/8} {\cal L}(\tau)\underline p & \leq & \epsilon_1 \dot \zeta(\tau)  + \eta^{-1} h(\tau)\left( \zeta e^{\eta^2/8} \phi_0' - q \right)  -\biggl(\dot q -\frac{3}{4} q\biggl) . 
\end{eqnarray}
For $\eta \in [1,\eta_1]$, where $\eta^{-1} < 1$, we have
\begin{eqnarray}\label{Lequivalent3}
\eta^{-1} e^{\eta^2/8} {\cal L}(\tau)\underline p & \leq & \epsilon_1 \dot \zeta(\tau)  + e^{-\tau/4} ( \zeta  + q)  -\biggl(\dot q -\frac{3}{4} q\biggl) .
\end{eqnarray}
To make this non-positive, we choose $\zeta$ to satisfy
\begin{eqnarray} 
 \epsilon_1 \dot \zeta(\tau)  \leq  \dot q -\frac{3}{4} q - e^{-\tau/4} ( \zeta + q) \label{zetareq1} 
 =   e^{-\tau/4} \zeta(\tau_0) -\frac{7}{4} q(\tau) - e^{-\tau/4} ( \zeta(\tau) + q(\tau)), %\no
\end{eqnarray}
where the last equalilty comes from (\ref{e5.3007}). Assuming $\dot \zeta < 0$, we have $\zeta(\tau) < \zeta(\tau_0)$, so a sufficient condition for (\ref{zetareq1}) to hold when $\tau \geq \tau_0$ is simply
\begin{eqnarray} \label{zetareq2}
 \epsilon_1 \dot \zeta(\tau) & \leq &  -3 q(\tau).
\end{eqnarray}
For $\eta$ near $0$, the dominant term in (\ref{Lequivalent2}) is $\eta^{-1} h(\tau)\left( \zeta e^{\eta^2/8} \phi_0' - q \right)$.  Define
\[
\epsilon_2 = \min_{\eta \in [0,1]} \phi_0'(\eta)e^{\eta^2/8} > 0.
\]
Therefore, if we can arrange that $\zeta(\tau) > q(\tau)/\epsilon_2$, then for $\eta \in [0,1]$, we have $\zeta e^{\eta^2/8} \phi_0' - q  \geq 0$, so
\[
 \eta^{-1} h(\tau)\left( \zeta e^{\eta^2/8} \phi_0' - q \right) \leq 0.
\]
In this case,
\begin{eqnarray}\label{Lequivalent4}
\eta^{-1} e^{\eta^2/8} {\cal L}(\tau)\underline p & \leq & \epsilon_1 \dot \zeta(\tau)  -\biggl(\dot q -\frac{3}{4} q\biggl) . 
\end{eqnarray}
which is non-positive for $\tau \geq \tau_0$, due to (\ref{zetareq1}). In summary, we will have ${\cal L}(\tau)\underline p \leq 0$ 
in the interval~$\eta \in [0,\eta_1]$ and $\tau \geq \tau_0$ if $\zeta$ satisfies (\ref{zetareq2}) and $\zeta(\tau) > q(\tau)/\epsilon_2$ for $\tau \geq \tau_0$. In view of this, we let $\zeta(\tau)$ have the form
\[
\zeta(\tau) = a_2 + a_3 e^{-(\tau - \tau_0)/4}.
\]
Thus, (\ref{zetareq2}) holds if
\[
- \frac{\epsilon_1 a_3}{4} e^{-(\tau - \tau_0)/4} \leq   - 3q = - 3e^{-(\tau - \tau_0)} - 4 e^{- \tau/4} (a_2 + a_3), \quad \tau \geq \tau_0.
\]
Hence it suffices that
\[
\frac{\epsilon_1 a_3}{4} \geq 3 + 4e^{-\tau_0/4} (a_2 + a_3)
\]
holds; this may be achieved with $a_2, a_3 > 0$ if $\tau_0$ is large enough. Then we may take $a_2$ large enough so that $\zeta(\tau) > q(\tau)/\epsilon_2$ also holds for $\tau \geq \tau_0$; this condition translates to:
\[
a_2 + a_3 e^{-(\tau - \tau_0)/4} \geq \frac{1}{\epsilon_2} \left(e^{-(\tau - \tau_0)} + \frac{4}{3}  e^{-\tau/4} (a_2 + a_3)\right), \quad \tau \geq \tau_0.
\]
This also is attainable with $a_2 > \frac{1}{\epsilon_2}$ and $a_3 > 0$ if $\tau_0$ is chosen large enough.  This completes the construction of the subsolution $\underline p(\tau,\eta)$ in (\ref{e5.3000}).

%%%%%%%%%%%%%%%%%%%end commentout

Let us come back to our subsolution $\underline z(\tau,\eta)$. From the strong maximum principle, 
we know that~$\underline z(\tau_0,\eta)>0$ and $\partial_\eta\underline z(\tau_0,0)>0$. Hence, there is $\lambda_0>0$ such that
\[
w(\tau_0,\eta)\geq\lambda_0\underline p(\tau_0,\eta),
\]
where $\underline p$ is given by (\ref{e5.3000}) with $\zeta$ and $q$ defined above, and we have for $\tau\geq\tau_0$:
$$
\underline w(\tau,\eta)\geq\lambda_0p(\tau,\eta).
$$
This, by (\ref{may204}), bounds the quantity $\underline I(\tau)$ uniformly from below, so that (\ref{may204}) holds with a constant~$c_0>0$
that depends on the initial condition $w_0$.

Therefore, just as in the study of the upper barrier, we obtain the uniform convergence of (possibly a subsequence of)  
$\underline w(\tau_n,\cdot)$ on the half-line $\eta\ge e^{-(1/2-\gamma)\tau}$ 
to a function~$\underline W_\infty(\eta)$ which satisfies
\begin{equation}\label{sep2512}
C^{-1}\eta e^{-\eta^2/4}\le \underline W_\infty(\eta)\le C \eta e^{-\eta^2/4},
\end{equation}
and such that
\begin{equation}\label{sep2510}
\lim_{n\to+\infty} e^{\eta^2/8}\vert\underline w(\tau_n,\eta)-\underline W_\infty(\eta)\vert=0, \quad \quad \eta > 0.
\end{equation}
\subsubsection*{Convergence of $w(\tau,\eta)$: proof of Lemma \ref{l3.1bis}}

Let $X$ be the space of bounded uniformly continuous functions $u(\eta)$ such that 
$e^{\eta^2/8}u(\eta)$ is bounded and uniformly continuous on $\R_+$. 
We deduce from the convergence of the upper and 
lower barriers for~$w(\tau,\eta)$ (and ensuing uniform bounds for $w$) that there
exists a sequence $\tau_n\to+\infty$ such that~$w(\tau_n,\cdot)$ itself
converges to a 
limit $W_\infty\in X$, such that $W_\infty\equiv0$ on $\R_-$, and $W_\infty(\eta)>0$ 
for all~$\eta>0$. Our next step is to 
bootstrap the convergence along a sub-sequence, and show that
the limit of $w(\tau,\eta)$ as $\tau\to+\infty$ exists in the space $X$. First, observe that the
above convergence implies that the shifted functions $
w_n(\tau,\eta)=w(\tau+\tau_n,\eta)$
converge in $X$, uniformly on compact time intervals,
as $n\to+\infty$ to the solution $w_\infty(\tau,\eta)$ of the linear problem
\begin{eqnarray}\label{sep2520}
&&(\partial_\tau+L)w_\infty=0,~~\eta>0,\\
&&w_\infty(\tau,0)=0,\nonumber\\
&&w_\infty(0,\eta)=W_\infty(\eta).\nonumber
\end{eqnarray}
In addition, there exists $\alpha_\infty>0$ such that $w_\infty(\tau,\eta)$ 
converges to 
$\bar\psi(\eta)=\alpha_\infty\eta e^{-\eta^2/4}$, 
in the topology of $X$ as $\tau\to+\infty$. Thus, for any $\eps>0$ we may 
%first choose $n$ sufficiently
%large so that
%\begin{equation}\label{sep2522}
%|w(\tau_n,\eta)-W_\infty(\eta)|\le \eps e^{-\eta^2/8},
%\end{equation}
%and then 
choose~$T_\eps$ large enough so that
\begin{equation}\label{sep2524}
|w_\infty(\tau,\eta)-\alpha_\infty\eta e^{-\eta^2/4}|\le \eps \eta e^{-\eta^2/8}
\hbox{ for all $\tau>T_\eps$, and $\eta>0$}.
\end{equation}
Given $T_\eps$ we can find $N_\eps$ sufficiently large so that
\begin{equation}\label{sep2526}
|w(T_\eps+\tau_n,\eta+e^{-(1/2-\gamma)T_\eps})-w_\infty(T_\eps,\eta)|\le \eps \eta e^{-\eta^2/8},
\hbox{ for all $n>N_\eps$.}
\end{equation}
In particular, we have
\begin{equation}\label{sep2528}
\alpha_\infty\eta e^{-\eta^2/4} -2\eps\eta e^{-\eta^2/8} 
\leq w(\tau_{N_\eps}+T_\eps,\eta+e^{-(1/2-\gamma)T_\eps})
\leq  \alpha_\infty \eta e^{-\eta^2/4}+2\eps\eta e^{-\eta^2/8}.
\end{equation}
We may now construct the upper and lower barriers for the function  $w(\tau+\tau_{N_\eps}+T_\eps,\eta+e^{-(1/2-\gamma)T_\eps})$, exactly as we have done before.
It follows, once again from Lemma~\ref{l3.2bis}   applied to these barriers  that any limit point $\phi_\infty$ 
of $w(\tau,\cdot)$ in $X$ as $\tau\to+\infty$ satisfies
\begin{equation}\label{sep2529}
(\alpha_\infty-C\eps)\eta e^{-\eta^2/4} 
\leq \phi_\infty(\eta)\leq  (\alpha_\infty+C\eps)\eta e^{-\eta^2/4}.
\end{equation}
As $\eps>0$ is arbitrary, we conclude that $w(\tau,\eta)$ converges in $X$ as $\tau\to+\infty$
to $\bar\psi(\eta) = \alpha_\infty \eta e^{-\eta^2/4}$. Taking into account Lemma~\ref{l3.2bis} once again, applied to the upper and lower barriers for $w(\tau,\eta)$
constructed starting from any time $\tau>0$, we have proved Lemma \ref{l3.1bis}, which implies Lemma \ref{sep29-lem02}.


\begin{thebibliography}{99}
\bibitem{ABBS} E. A\"id\'ekon, J. Berestycki, \'E. Brunet, Z. Shi, {\em Branching Brownian motion seen from its tip}, Probab. Theory Relat. Fields {\bf 157} (2013), pp. 405-451.

\bibitem{ABK1} L.-P. Arguin, A. Bovier, and N. Kistler, {\em Poissonian statistics in the extremal process of branching Brownian motion}. Ann. Appl. Probab. {\bf 22} (2012), pp. 1693-1711.


\bibitem{ABK2} L.-P. Arguin, A. Bovier, and N. Kistler, {\em The extremal process of branching Brownian motion}. Probab. Theory Relat. Fields {\bf 157} (2013) pp. 535-574.

\bibitem{BH} H. Berestycki, F. Hamel, {Front propagation in periodic excitable media}, Comm. Pure Appl. Math. {\bf 55} (2002), 949--1032.


\bibitem{Bramson1} 
M.D. Bramson, {Maximal displacement of branching Brownian motion}, 
Comm. Pure Appl. Math. {\bf 31}, 1978, 531--581.

\bibitem{Bramson2} 
M.D. Bramson, {Convergence of solutions of the Kolmogorov equation
  to travelling waves}, Mem. Amer. Math. Soc. {\bf 44}, 1983. 



\bibitem{BD2} E. Brunet and B. Derrida.  {\em A branching random walk seen from the tip}, Journal of Statistical Physics. {\bf 143} (2011), pp. 420-446.


\bibitem{BD1} E. Brunet and B. Derrida. {\em Statistics at the tip of a branching random walk and the delay of traveling waves}. Eur. Phys. Lett. 87, 60010 (2009).


\bibitem{FZ} M. Fang and O. Zeitouni, {Slowdown for time inhomogeneous branching 
Brownian motion},  J.~Stat. Phys. {\bf 149}, 2012, 1--9.

\bibitem{Fisher} R.A. Fisher, {The wave of advance of advantageous genes}, Ann. Eugenics {\bf 7}, 1937, 353--369.

\bibitem{HNRR1} F. Hamel, J. Nolen, J.-M. Roquejoffre and L. Ryzhik,
  {A short proof of the logarithmic Bramson correction in Fisher-KPP
    equations}, Netw. Het. Media {\bf 8}, 2013, 275--289.

\bibitem{HNRR2} F. Hamel, J. Nolen, J.-M. Roquejoffre, and L. Ryzhik, {The 
logarithmic time delay of KPP fronts in a periodic medium,} J. Europ. Math. 
Soc. {\bf 18}, 2016, 465--505.

\bibitem{HR} F. Hamel, L. Roques, {Uniqueness and stability properties 
of monostable pulsating fronts}, J.~Europ. Math. Soc. {\bf 13}, 2011, 345--390. 

\bibitem{Henderson} C. Henderson, Population stabilization in branching 
Brownian motion with absorption, to appear in CMS, 2015.

\bibitem{KPP} A.N. Kolmogorov, I.G. Petrovskii and N.S. Piskunov, {\'Etude de 
l'\'equation de la diffusion avec croissance de la quantit\'e de mati\`ere et son 
application \`a un probl\`eme biologique}, Bull. Univ. \'Etat Moscou, S\'er. 
Inter.~A {\bf 1}, 1937, 1--26.

\bibitem{LS} S.P. Lalley and T. Sellke, { A conditional limit theorem for the frontier of a branching Brownian motion}. Annals of Probability, {\bf 15},
1987, 1052--1061.


\bibitem{Lau} K.-S. Lau, {On the nonlinear diffusion equation of Kolmogorov, Petrovskii and Piskunov}, J.~Diff. Eqs. {\bf 59}, 1985, 44-70.

\bibitem{MZ} P. Maillard, O. Zeitouni, Slowdown in branching Brownian 
motion with inhomogeneous variance, to appear in Ann. IHP, Prob. Stat.
  
\bibitem{McK}   H.P. McKean, {Application of Brownian motion 
to the equation of Kolmogorov-Petrovskii-Piskunov}, Comm. Pure Appl. Math. {\bf 28} 
1975, 323--331.


\bibitem{NRR1} J. Nolen, J.-M. Roquejoffre and L. Ryzhik, 
{Power-like delay in time inhomogeneous Fisher-KPP equations,} Comm. Partial Diff. 
Equations, 
{\bf 40}, 2015, 475--505 

\bibitem{NRR2} J. Nolen, J.-M. Roquejoffre and L. Ryzhik, Sharp large-time asymptotics
in the Fisher-KPP equation, forthcoming.

\bibitem{Roberts} M. Roberts, 
{A simple path to asymptotics for
the frontier of a branching Brownian motion}, Ann. Prob. {\bf 41}, 2013, 
3518--3541.

\bibitem{Roq97} J.-M. Roquejoffre, Eventual monotonicity and convergence to 
travelling fronts for the solutions of parabolic equations in cylinders, 
Ann. Inst. H. Poincar\'e Anal. Non Lin\'eaire {\bf 14}, 1997, 
499--552. 

\end{thebibliography}
\end{document}